%% file: MatrixCEN.tex
\documentclass{article} 
\usepackage{nips12submit_e,times}

\usepackage[dvips]{graphics}
\DeclareGraphicsExtensions{.eps.gz,.eps,.epsi.gz,.epsi,.ps,.ps.gz}
\usepackage{epsfig}
\usepackage{color}
\graphicspath{{fig/}} 

\usepackage{latexsym}
\usepackage{graphicx}
\usepackage{amsfonts}
\usepackage{amssymb}
\usepackage{mathrsfs}
\usepackage{verbatim}

\input{04-pre-am-short.tex}

\usepackage{amsmath}
\usepackage{latexsym}
\usepackage{mathrsfs}
\usepackage{amsfonts}
\usepackage{amssymb}

\newtheorem{theorem}{Theorem}
\newtheorem{lemma}{Lemma}
\newtheorem{proposition}{Proposition}

\newtheorem{algorithm}{Algorithm}

\newcommand{\chzblue}[1]{{{#1}}}

\def\argmin{\mathop{\rm arg\, min}}
\def\real{\mathop{{\rm I}\kern-.2em\hbox{\rm R}}\nolimits}

\title{Calibrated Elastic Regularization in Matrix Completion}

\author{
Tingni Sun\\
Statistics Department, The Wharton School\\
University of Pennsylvania\\
Philadelphia, Pennsylvania 19104\\
\texttt{tingni@wharton.upenn.edu} \\
\And
Cun-Hui Zhang\\
Department of Statistics and Biostatistics\\
Rutgers University\\
Piscataway, New Jersey 08854\\
\texttt{czhang@stat.rutgers.edu} \\
}

\nipsfinalcopy 

\begin{document}

\maketitle

\begin{abstract}
This paper concerns the problem of matrix completion, which is to estimate a matrix from observations in a small subset of indices. We propose a calibrated spectrum elastic net method with a sum of the nuclear and Frobenius penalties and develop an iterative algorithm to solve the convex minimization problem. The iterative algorithm alternates between imputing the missing entries in the incomplete matrix by the current guess and estimating the matrix by a scaled soft-thresholding singular value decomposition of the imputed matrix until the resulting matrix converges. A calibration step follows to correct the bias caused by the Frobenius penalty. Under proper coherence conditions and for suitable penalties levels,
we prove that the proposed estimator achieves an error bound of nearly optimal order and in proportion to the noise level.
This provides a unified analysis of the noisy and noiseless matrix completion problems.
Simulation results are presented to compare our proposal with previous ones.
\end{abstract}

\section{Introduction}
Let $\Theta\in\real^{d_1\times d_2}$ be a matrix of interest and $\Omega^*=\{1,\ldots,d_1\}\times \{1,\ldots,d_2\}$.
Suppose we observe vectors $(\omega_i,y_i)$,
\bel{model-iid}
y_i = \Theta_{\omega_i} + \veps_i,\quad i=1,\dots,n,
\eel
where $\omega_i\in \Omega^*$ and $\veps_i$ are random errors. We are interested in estimating $\Theta$ when
$n$ is a small fraction of $d_1d_2$.
A well-known application of matrix completion is the Netflix problem where $y_i$ is the rating of movie $b_j$ by user $a_i$ for $\omega=(a_i,b_j)\in\Omega^*$ \cite{Netflix07}. In such applications, the proportion of the observed entries is typically very small, so that the estimation or recovery of $\Theta$ is impossible without a structure assumption on $\Theta$. In this paper, we assume that $\Theta$ is of low rank.

A focus of recent studies of matrix completion has been on a simpler formulation, also known as exact recovery, where the observations are assumed to be uncorrupted, i.e. $\veps_{i}=0$. A direct approach is to minimize $\rank(M)$ subject to $M_{\omega_i}=y_i$. 
\chzblue{An iterative algorithm was proposed in \cite{KesMO09-exact} to project a trimmed SVD of the 
incomplete data matrix to the space of matrices of a fixed rank $r$.  
The }nuclear norm was proposed as a surrogate for the rank, leading to the following convex
minimization problem in a linear space  \cite{CandesR09}:
\bes
\Thetahat^{(\text{CR})} = \argmin_{M}\Big\{\|M\|_{(N)}: M_{\omega_i}=y_i \,\forall\, i\le n\Big\}.
\ees
We denote the nuclear norm by $\| \cdot\|_{(N)}$ here and throughout this paper.
This procedure, analyzed in \cite{CandesR09, CandesT09, Gross09, Retch09} among others,
is parallel to the replacement of the $\ell_0$ penalty by the $\ell_1$ penalty in solving
the sparse recovery problem in a linear space. 

In this paper, we focus on the problem of matrix completion with noisy observations (\ref{model-iid}) and take 
the exact recovery as a special case. Since the exact constraint \chzblue{is no longer appropriate in the 
presence of noise, penalized squared error $\sum_{i=1}^n(M_{\omega_i}-y_i)^2$} is considered.
By reformulating the problem in Lagrange form, \cite{MazHT10} proposed the spectrum Lasso
\bel{matrix-lasso}
\Thetahat^{(\text{MHT})} = \argmin_{M}\Big\{
\sum_{i=1}^n M_{\omega_i}^2/2 - \sum_{i=1}^n y_i M_{\omega_i} + \lambda \| M\|_{(N)}\Big\},
\eel
along with an iterative convex minimization algorithm.
However, (\ref{matrix-lasso}) is difficult to analyze when the sample fraction $\pi_0=n/(d_1d_2)$ is small,
due to the ill-posedness of the quadratic term $\sum_{i=1}^n M_{\omega_i}^2$.
This has led to two alternatives in \cite{KolTL10} and \cite{NegW10}.
While \cite{NegW10} proposed to minimize (\ref{matrix-lasso}) under an additional $\ell_\infty$ constraint on $M$,
\cite{KolTL10} modified (\ref{matrix-lasso}) by replacing the quadratic term $\sum_{i=1}^n M_{\omega_i}^2$ with
$\pi_0\|M\|_{(F)}^2$.
Both \cite{KolTL10, NegW10} provided nearly optimal error bounds when the noise level is of no smaller order
than the $\ell_\infty$ norm of the target matrix $\Theta$, but not of smaller order,
especially not for exact recovery.
In a different approach, \cite{KesMO09} proposed a non-convex recursive algorithm and
provided error bounds in proportion to the noise level. However, the procedure requires the
knowledge of the rank $r$ of the unknown $\Theta$ and the error bound is optimal only
when $d_1$ and $d_2$ are of the same order.

Our goal is to develop an algorithm for matrix completion that can be as easily computed as the spectrum Lasso
(\ref{matrix-lasso}) and enjoys a nearly optimal error bound proportional to the noise level to continuously cover
both the noisy and noiseless cases.
We propose to use an elastic penalty, a linear combination of the nuclear and Frobenius norms, which leads to the estimator
\bel{Thetatil}
\Thetatil = \argmin_M\Big\{ \sum_{i=1}^n M_{\omega_i}^2/2 - \sum_{i=1}^n y_i M_{\omega_i}
+ \lam_1 \| M\|_{(N)} + (\lam_2/2)\|M\|^2_{(F)}\Big\},
\eel
where $\|\cdot\|_{(N)}$ and $\| \cdot\|_{(F)}$ are the nuclear and Frobenius norms, respectively.
We call (\ref{Thetatil}) spectrum elastic net (E-net) since it is parallel to the E-net in linear regression,
the least squares estimator with a sum of the $\ell_1$ and $\ell_2$ penalties, introduced in \cite{ZouH05}.
Here the nuclear penalty provides the sparsity in the spectrum, while the Frobenius penalty regularizes
the inversion of the quadratic term.
Meanwhile, since the Frobenius penalty roughly shrinks the estimator by a factor $\pi_0/(\pi_0+\lam_2)$,
we correct this bias by a calibration step,
\bel{cali-elastic}
\Thetahat = (1+\lam_2/\pi_0) \Thetatil.
\eel
We call this estimator calibrated spectrum E-net.

Motivated by \cite{MazHT10}, we develop an EM algorithm to solve (\ref{Thetatil}) for matrix completion.
The algorithm iteratively replaces the missing entries with those obtained from a scaled soft-thresholding singular value decomposition (SVD) until the resulting matrix converges. This EM algorithm is guaranteed to converge to
the solution of (\ref{Thetatil}).

Under proper coherence conditions, we prove that for suitable penalty levels $\lam_1$ and $\lam_2$,
the calibrated spectrum E-net (\ref{cali-elastic}) achieves a desired error bound in the Frobenius norm.
Our error bound is of nearly optimal order and in proportion to the noise level.
This provides a sharper result than those of \cite{KolTL10,NegW10} when the
noise level is of smaller order than the $\ell_\infty$ norm of $\Theta$, and than that of
\cite{KesMO09} when $d_2/d_1$ is large. Our simulation results support the use of the calibrated spectrum E-net.
They illustrate that (\ref{cali-elastic}) performs comparably to (\ref{matrix-lasso}) and outperforms
the modified method of \cite{KolTL10}.

Our analysis of the calibrated spectrum E-net uses an inequality similar to a duel certificate bound in \cite{CandesT09}.
The bound in \cite{CandesT09} requires sample size $n\asymp \min\{(r\log d)^2, r(\log d)^6\}d\log d$,
where $d = d_1+d_2$. We use the method of moments to remove a $\log d$ factor in the first component of
their sample size  requirement. This leads to a sample size requirement of $n \asymp r^2d\log d$, with an extra $r$
in comparison to the \chzblue{ideal} $n \asymp rd\log d$. Since the extra $r$ does not appear in our error bound,
its appearance in the sample size requirement seems to be a technicality.

The rest of the paper is organized as follows. 
In Section 2, we describe an iterative algorithm for the computation of the spectrum E-net 
\chzblue{and study its convergence}. 
In Section 3, we derive error bounds for the calibrated spectrum E-net. Some simulation results are presented in Section 4. Section 5 provides the proof of our main result.

We use the following notation throughout this paper.
For matrices $M\in\R^{d_1\times d_2}$, $\|M\|_{(N)}$ is the nuclear norm (the sum of all singular values of $M$), $\|M\|_{(S)}$ is the spectrum norm (the largest singular value), $\|M\|_{(F)}$ is the Frobenius norm 
\chzblue{(the $\ell_2$ norm of vectorized $M$)}, 
and \chzblue{$\|M\|_\infty=\max_{jk}|M_{jk}|$}. 
Linear mappings from $\R^{d_1\times d_2}$ to $\R^{d_1\times d_2}$ are denoted by the calligraphic letters. For a linear mapping $\calQ$, the operator norm is $\|\calQ\|_{(op)}=\sup_{\|M\|_{(F)}=1}\|\calQ M\|_{(F)}$.
We equip $\R^{d_1\times d_2}$ with the inner product $\langle M_1,M_2\rangle = \trace(M_1^\top M_2)$ so that $\langle M,M\rangle=\|M\|_{(F)}^2$. 
\chzblue{For projections $\calP$, $\calP^\perp = \calI - \calP$ with $\calI$ being the identity}. 
We denote by $E_\omega$ the unit matrix with 1 at $\omega\in\{1,\ldots,d_1\}\times \{1,\ldots,d_2\}$,
and by $\calP_\omega$ the projection to $E_\omega$: $M \to M_\omega E_\omega = \langle E_\omega,M\rangle E_\omega$.

\section{An algorithm for spectrum elastic regularization}
We first present a lemma for the M-step of our iterative algorithm.
\begin{lemma}\label{lem-threshold}
Suppose the matrix $Z$ has rank $r$. The solution to the optimization problem
\bes
\argmin_{Z} \Big\{\|Z-W\|_{(F)}^2/2+\lam_1\|Z\|_{(N)}+\lam_2\|Z\|_{(F)}^2/2\Big\}
\ees
is given by $S(W;\lam_1,\lam_2)=UD_{\lam_1,\lam_2}V'$ with $D_{\lam_1,\lam_2}=\diag\{(d_1-\lam_1)_+,\dots,(d_r-\lam_1)_+\}/(1+\lam_2)$, where $UDV'$ is the SVD of $W$, $D=\diag\{d_1,\dots,d_r\}$ and $t_+=max(t,0)$.
\end{lemma}
The minimization problem in Lemma \ref{lem-threshold} is solved by a scaled soft-thresholding SVD.
This is parallel to Lemma 1 in \cite{MazHT10} and justified by Remark 1 there. We use Lemma \ref{lem-threshold}
to solve the M-step of the EM algorithm for the spectrum E-net (\ref{Thetatil}).

We still need an E-step to impute a complete matrix given the observed data $\{y_i, \omega_i: i=1,\dots,n\}$.
Since $\omega_i$ are allowed to have ties, we need the following notation.
Let $m_\omega=\#\{i:\omega_i=\omega,i\le n\}$ be the multiplicity of observations at $\omega\in \Omega^*$
and $m^*=\max_{\omega}m_\omega$ be the maximum multiplicity. Suppose that the complete data is
composed of $m_*$ observations at each $\omega$ for a certain integer $m_*$. Let $\Ybar^{(\text{com})}_\omega$
be the sample mean of the complete data at $\omega$ and $\Ybar^{(\text{com})}$ be the matrix with
components $\Ybar^{(\text{com})}_\omega$.  If the complete data are available, (\ref{Thetatil}) is equivalent to
\bes
\argmin_M\Big\{ (m_*/2) \|\Ybar^{(\text{com})}-M\|_{(F)}^2 + \lam_1 \| M\|_{(N)} + (\lam_2/2)\|M\|^2_{(F)}\Big\}.
\ees
Let $\Ybar^{(\text{\rm obs})}_\omega = m_\omega^{-1}\sum_{\omega_i=\omega} y_i$ be the sample mean of
the observations at $\omega$ and $\Ybar^{(\text{\rm obs})} = (\Ybar^{(\text{\rm obs})}_\omega)_{d_1\times d_2}$.
In the white noise model, the conditional expectation of $\Ybar^{(\text{com})}_\omega$ given $\Ybar^{(\text{\rm obs})}$ is
$(m_\omega/m_*)\Ybar^{(\text{\rm obs})}_\omega + (1-m_\omega/m_*)\Theta_\omega$ for $m_\omega\le m_*$.
This leads to a generalized E-step:
\bel{y-exp}
\Ybar^{(\text{\rm imp})} = (\Ybar^{(\text{\rm imp})}_\omega)_{d_1\times d_2},\
\Ybar^{(\text{\rm imp})}_\omega = \min\{1,(m_\omega/m_*)\}\Ybar^{(\text{\rm obs})}_\omega + (1-m_\omega/m_*)_+Z^{(\text{\rm old})}_\omega,
\eel
where $Z^{(\text{\rm old})}$ is the estimation of $\Theta$ in the previous iteration. This is a genuine E-step when
$m_*=m^*$ but also allows a smaller $m_*$ to reduce the proportion of missing data.

We now present the EM-algorithm for the computation of the spectrum E-net $\Thetatil$ in
(\ref{Thetatil}).

\begin{algorithm}\label{alg}
Initialize with $Z^{(0)}$ and $k=0$. Repeat the following steps:
\begin{itemize}
\item E-step: Compute $\Ybar^{(\text{\rm imp})}$ in (\ref{y-exp}) with $Z^{(\text{\rm old})} = Z^{(k)}$ and assign $k\leftarrow k+1$,
\item M-step: Compute $Z^{(k)}=S(\Ybar^{(\text{\rm imp})};\lam_1/m_*,\lam_2/m_*)$,
\end{itemize}
until $\|Z^{(k)}-Z^{(k-1)}\|_{(F)}^2/\|Z^{(k)}\|_{(F)}^2\le \eps$. Then, return $Z^{(k)}$.
\end{algorithm}

The following theorem states the convergence of Algorithm \ref{alg}.
\begin{theorem}\label{thm-converge} As $k\to\infty$, $Z^{(k)}$ converges to a limit $Z^{(\infty)}$ as a
function of the data and $(\lam_1,\lam_2,m_*)$, and $Z^{(\infty)}=\Thetatil$ for $m_*\ge m^*$.
\end{theorem}

Theorem \ref{thm-converge} is a variation of a parallel result in \cite{MazHT10} and follows from the same
proof there. As \cite{MazHT10} pointed out, a main advantage of Algorithm \ref{alg} is the speed of each iteration.
When the maximum multiplicity $m^*$ is small, we simply use $Z^{(0)}=\Ybar^{(\text{\rm obs})}$ and $m_*=m^*$;
Otherwise, we may first run the EM-algorithm for an $m_*<m^*$ and use the output as the initialization $Z^{(0)}$
for a second run of the EM-algorithm with $m_*=m^*$.

\section{Analysis of estimation accuracy}
In this section, we derive error bounds for the calibrated spectrum E-net.
We need the following notation. Let $r=\rank(\Theta)$, $UDV^\top$ be the SVD of $\Theta$, and $s_1\ge \ldots \ge s_r$
be the nonzero singular values of $\Theta$.
Let $T$ be the tangent space with respect to ${U}{V}^\top$, the space of all matrices of the form
$UU^\top M_1+M_2VV^\top$. The orthogonal projection to $T$ is given by
\bel{matrix-completion-P_T}
\calP_TM = UU^\top M + MVV^\top  - UU^\top MVV^\top .
\eel

\begin{theorem}\label{thm-oracle}
Let ${\xi}=1+\lam_2/\pi_0$ and $\calH = \sum_{i=1}^n\calP_{\omega_i}$. Define
\bes
\calR &=& (\calH - \pi_0)\calP_T/(\pi_0+\lam_2),\\
\Deltabar &=&  \calR(\lam_2 \Theta + \lam_1 {U}{V}^\top),\\
\calQ  &=& \calI -  \calH(\calP_T\calH\calP_T+\lam_2\calP_T)^{-1}\calP_T.
\ees
Let $\veps = \sum_{i=1}^n \veps_i E_{\omega_i}$. Suppose
\begin{align}
&\|\calP_T\calR\|_{(op)}\le 1/2, \quad s_r\ge 5\lam_1/\lam_2, \label{event1}\\
&\|\calP_T\Deltabar\|_{(F)}\le\sqrt{r}\lam_1/8,\quad
\|\Deltabar - \calR(\calP_T\calR+\calP_T)^{-1}\calP_T\Deltabar\big\|_{(S)} \le \lam_1/4, \label{event2}\\
&\|\calP_T\veps\|_{(F)}\le \sqrt{r}\lam_1/8, \quad
 \|\calQ\veps\|_{(S)}\le 3\lam_1/4, \quad
 \|\calP_T^\perp\veps\|_{(S)}\le \lam_1.  \label{event3}
\end{align}
Then the calibrate spectrum E-net (\ref{cali-elastic}) satisfies
\bel{oracle-bound}
\|\Thetahat - \Theta\|_{(F)} \le 2\sqrt{r}\lam_1/\pi_0.
\eel
\end{theorem}


The proof of Theorem \ref{thm-oracle} is provided in Section 5.
\chzblue{When $\omega_i$ are random entries in $\Omega^*$, $E\calH = \pi_0\calI$, so that}
(\ref{event2}) and the first inequality of
(\ref{event1}) \chzblue{are expected to} hold under proper conditions.
Since the rank of $\calP_T\veps$ is no greater than $2r$,
(\ref{event3}) essentially requires $\|\veps\|_{(S)}\asymp \lam_1$.
Our analysis allows $\lam_2$ to lie in a certain range $[\lam_*,\lam^*]$, and $\lam^*/\lam_*$
is large under proper conditions. 
\chzblue{Still, the choice of $\lam_2$ is constrained by (\ref{event1}) and (\ref{event2}) 
since $\Deltabar$ is linear in $\lam_2$}. 
When $\lam_2/\pi_0$ diverges to infinity, the calibrated spectrum E-net (\ref{cali-elastic}) becomes the modified
spectrum Lasso of \chzblue{\cite{KolTL10}}.

Theorem \ref{thm-oracle} provides sufficient conditions on the target matrix and the noise for achieving a certain level
of estimation error. Intuitively, these conditions on the target matrix $\Theta$ must imply a certain level of coherence
(or flatness) of the unknown matrix since it is impossible to distinguish the unknown from zero when the observations
are completely outside its support.
In \cite{CandesR09, CandesT09, Gross09, Retch09}, coherence conditions are imposed on
\bel{coherence}
\mu_0 = \max\{ (d_1/r)\|UU^\top\|_\infty, (d_2/r)\|VV^\top\|_\infty\},\quad \mu_1 = \sqrt{d_1d_2/r}\|UV^\top\|_\infty,
\eel
where $U$ and $V$ are matrices of singular vectors of $\Theta$.
\cite{NegW10} considered a more general notation of spikiness of a matrix $M$, defined as the ratio between
the $\ell_\infty$ and dimension-normalized $\ell_2$ norms,
\bel{spikiness}
\alpha_{sp}(M)=\|M\|_\infty\sqrt{d_1d_2}/\|M\|_{(F)}.
\eel

Suppose in the rest of the section that $\omega_i$ are iid points uniformly distributed in $\Omega^*$
and $\veps_i$ are iid $N(0,\sigma^2)$ variables independent of $\{\omega_i\}$.
The following theorem asserts that under certain coherence conditions on the matrices $\Theta$,
$UU^\top$, $VV^\top$ and $UV^\top$, all conditions of Theorem \ref{thm-oracle} hold with large probability
when the sample size $n$ is of the order $r^2 d\log d$.

\begin{theorem}\label{thm-noise} Let $d=d_1+d_2$. Consider $\lam_1$ and $\lam_2$ satisfying
\bel{lams}
\lam_1 = \sigma\sqrt{8\pi_0 d\log d},\quad
1\le \frac{\lam_2\|\Theta\|_{(F)}}{\lam_1\{n/(d\log d)\}^{1/4}} \le 2.
\eel
Then, there exists a constant $C$ such that
\bel{samplesize}
n \ge C\max\Big\{\mu_0^2r^2 d\log d,(\mu_1+r)\mu_1r d\log d, (\alpha_{sp}^{4/3}\vee \kappa_*^4)
r^2d\log d\Big\}
\eel
implies
\bes
\|\Thetahat - \Theta\|_{(F)}^2/(d_1d_2) \le 32(\sigma^2 rd\log d)/n
\ees
with probability at least $1-1/d^2$, where $\mu_0$ and $\mu_1$ are the coherence constants in (\ref{coherence}), $\alpha_{sp}=\alpha_{sp}(\Theta)$ is the spikiness of $\Theta$ and $\kappa_* = \|\Theta\|_{(F)}/(r^{1/2}s_r)$.
\end{theorem}

We require the knowledge of noise level $\sigma$ to determine the penalty level that is usually 
considered as tuning parameter in practice.
The Frobenius norm $\|\Theta\|_{(F)}$ in (\ref{lams}) can be replaced by
an estimate of the same magnitude in Theorem \ref{thm-noise}. In our simulation experiment, we use
$\lam_2 = \lam_1\{n/(d\log d)\}^{1/4}/\Fhat$ with $\Fhat = (\sum_{i=1}^n y_i^2/\pi_0)^{1/2}$.
\chzblue{The Chebyshev inequality} provides $\Fhat/\|\Theta\|_{(F)}\to 1$ when $\alpha_{sp}=O(1)$ and
\chzblue{$\sigma^2 \ll \|\Theta\|_\infty^2$}.


A key element in our analysis is to find a probabilistic bound for the second inequality of (\ref{event2}),
or equivalently an upper bound for
\bel{prob1}
P\big\{ \|\calR(\calP_T\calR+\calP_T)^{-1}(\lam_2 \Theta + \lam_1 {U}{V}^\top)\|_{(S)} > \lam_1/4 \big\}.
\eel
\chzblue{This guarantees the existence of a primal dual certificate for the spectrum E-net penalty \cite{ZhangZ12}}. 
For $\lam_2=0$, a similar inequality was proved in \cite{CandesT09},
where the sample size requirement is $n\ge C_0\min\{\mu^2 r^2(\log d)^2d, \mu^2 r(\log d)^6d\}$
for a certain coherence factor $\mu$.
We remove a log factor in the first bound, resulting in the sample size requirement in (\ref{samplesize}), 
\chzblue{which is optimal when $r=O(1)$}. 
For exact recovery in the noiseless case, the sample size 
$n\asymp rd(\log d)^2$ \chzblue{is sufficient if a golfing scheme is used to construct} an approximate 
dual certificate \cite{Gross09,Retch09}. 
We use the following lemma to bound (\ref{prob1}).

\vspace{-0.05in}
\begin{lemma}\label{lem-specbound}
Let $\calH = \sum_{i=1}^n\calP_{\omega_i}$ where $\omega_i$ are iid points uniformly distributed in $\Omega^*$.
Let $\calR = (\calH - \pi_0)\calP_T/(\pi_0+\lam_2)$ and ${\xi}=1+\lam_2/\pi_0$. Let $M$ be a deterministic matrix.
Then, there exists a numerical constant $C$ such that, for all $k\ge 1$ and $m\ge 1$,
\bel{V-stat}
{\xi}^{2km}E\|\calR^k M\|_{(S)}^{2m}
\le \Big\{ C\mu_0^2r^2 d km/n \Big\}^{km}\Big(\mu_0^{-2}(\sqrt{d_1d_2}/r)\|M\|_\infty\Big)^{2m}.
\eel
\end{lemma}

\vspace{-0.05in}
We use a different graphical approach than those in \cite{CandesT09} to bound
$E\,\trace(\{(\calR^kM)^\top(\calR^kM)\}^m)$ in the proof of Lemma \ref{lem-specbound}.
The rest of the proof of Theorem \ref{thm-noise} can be outlined as follows.
Assume that all coherence factors are $O(1)$.
Let $M = \lam_2\Theta+\lam_1 UV^\top$ and write
$\calR(\calP_T\calR+\calP_T)^{-1}M = \calR M - \calR^2M + \cdots + (-1)^{k^*-1}\calR^{k^*}M+\Rem$.
By (\ref{V-stat}) with $km\asymp \log d$ for $k\ge 2$ and an even simpler bound for $k=1$ and $\Rem$,
(\ref{prob1}) holds when $(\sqrt{d_1d_2}/r)\|M\|_\infty\asymp \lam_1\eta$, where $\eta \asymp r^2 d(\log d)/n$.
Since $\alpha_{sp} + \mu_1+\|\Theta\|_{(F)}^2/(rs_r^2)=O(1)$, this is equivalent to $\eta(s_r\lam_2/\lam_1+1) \lesssim 1$.
Finally, we use matrix exponential inequalities \cite{Oliveira10,Tropp11b} to verify other conditions of
Theorem~\ref{thm-oracle}.
We omit technical details of the proof of Lemma \ref{lem-specbound} and
Theorem \ref{thm-noise}. We would like to point out that if the $r^2$ in (\ref{V-stat}) can be replaced by
$r(\log d)^\gamma$, \chzblue{e.g. $\gamma = 5$ in view of \cite{CandesT09}, the rest of the proof of 
Theorem \ref{thm-noise} is intact with $\eta \asymp rd(\log d)^{1+\gamma}/n$ and a proper adjustment of 
$\lam_2$ in (\ref{lams})}. 

Compared with \cite{KolTL10} and \cite{NegW10}, the main advantage of Theorem \ref{thm-noise} is the
proportionality of its error bound to the noise level.
In \cite{KolTL10}, the quadratic term $\sum_{i=1}^n M_{\omega_i}^2$ in (\ref{matrix-lasso}) is replaced by its expectation $\pi_0\|M\|_{(F)}^2$ and the resulting minimizer is proved to satisfy
\bel{KLT-bd}
\|\Thetahat^{(\text{KLT})} - \Theta\|_{(F)}^2/(d_1d_2) \le C \max(\sigma^2,\|\Theta\|_\infty^2) rd(\log d)/n
\eel
with large probability, where $C$ is a numerical constant.
This error bound achieves the squared error rate $\sigma^2 rd(\log d)/n$ as in Theorem \ref{thm-noise}
when the noise level $\sigma$ is of no smaller order than $\|\Theta\|_\infty$, but not of smaller order.
In particular, (\ref{KLT-bd}) does not imply exact recovery when $\sigma=0$.
In Theorem~\ref{thm-noise}, the error bound converges to zero as the noise level diminishes,
implying exact recovery in the noiseless case.
In \cite{NegW10}, a constrained spectrum Lasso was proposed
that minimizes (\ref{matrix-lasso}) subject to $\|M\|_\infty \le \alpha^*/\sqrt{d_1d_2}$.
For $\|\Theta\|_{(F)}\le 1$ and $\alpha_{sp}(\Theta)\le \alpha^*$, \cite{NegW10} proved
\bel{NW-bd}
\|\Thetahat^{(\text{NW})} - \Theta\|_{(F)}^2 \le C \max(d_1d_2 \sigma^2, 1) (\alpha^*)^2 rd(\log d)/n
\eel
with large probability. Scale change from the above error bound yields
\bes
\|\Thetahat^{(\text{NW})} - \Theta\|_{(F)}^2/(d_1d_2) \le C \max\{\sigma^2, \|\Theta\|_{(F)}^2/(d_1d_2)\}
(\alpha^*)^2 rd(\log d)/n.
\ees
Since $\alpha^*\ge 1$ and $\alpha^*\|\Theta\|_{(F)}/\sqrt{d_1d_2}\ge \|\Theta\|_\infty$, the right-hand side
of (\ref{NW-bd}) is of no smaller order than that of (\ref{KLT-bd}). We shall point out that (\ref{KLT-bd}) and (\ref{NW-bd}) only
require sample size $n\asymp rd\log d$. 
In addition, \cite{NegW10} allows more practical weighted sampling models.

Compared with \cite{KesMO09}, the main advantage of Theorem \ref{thm-noise} is the independence
of its sample size requirement on the aspect ratio $d_2/d_1$, where $d_2\ge d_1$ is assumed without
loss of generality by symmetry. The error bound in \cite{KesMO09} implies
\bel{KMO-bd}
\|\Thetahat^{(\text{KMO})} - \Theta\|_{(F)}^2/(d_1d_2)
\le C_0 (s_1/s_r)^4 \sigma^2 rd(\log d)/n
\eel
for sample size $n\ge C^*_1 rd\log d + C^*_2 r^2 d \sqrt{d_2/d_1}$,
where $\{C^*_1, C^*_2\}$ are constants depending on the same set of coherence factors as in
(\ref{samplesize}) and $s_1>\cdots>s_r$ are the singular values of $\Theta$.
Therefore, Theorem \ref{thm-noise} effectively replaces the root aspect ratio $\sqrt{d_2/d_1}$ in
the sample size requirement of (\ref{KMO-bd}) with a log factor, and removes the coherence factor
$(s_1/s_r)^4$ on the right-hand side of (\ref{KMO-bd}). 
We note that $s_1/s_r$ is a larger coherence
factor than $\|\Theta\|_{(F)}/(r^{1/2}s_r)$ in the sample size requirement in Theorem \ref{thm-noise}.
\chzblue{The root aspect ratio can be removed from the sample size requirement for (\ref{KMO-bd}) 
if $\Theta$ can be divided into square blocks uniformly satisfying the coherence conditions}.

\section{Simulation study}
This experiment has the same setting as in Section 9 of \cite{MazHT10}.
We provide the description of the simulation settings in our notation as follows:
The target matrix is $\Theta=UV^\top$, where $U_{d_1\times r}$ and $V_{d_2\times r}$ are random matrices
with independent standard normal entries.
The sampling points $\omega_i$ have no tie and $\Omega=\{\omega_i:i=1,\dots,n\}$ is a uniformly
distributed random subset of $\{1,\dots,d_1\}\times\{1,\dots,d_2\}$, where $n$ is fixed.
The errors $\veps$ are iid $N(0,\sigma^2)$ variables.
Thus, the observed matrix is $Y=\calP_\Omega(\Theta+\veps)$ with
$\calP_\Omega = \calH = \sum_{i=1}^n\calP_{\omega_i}$ being a projection.
The signal to noise ratio (SNR) is defined as $\text{SNR}=\sqrt{r}/\sigma$.

\begin{figure}[htbp]
\includegraphics[width=.8\textwidth]{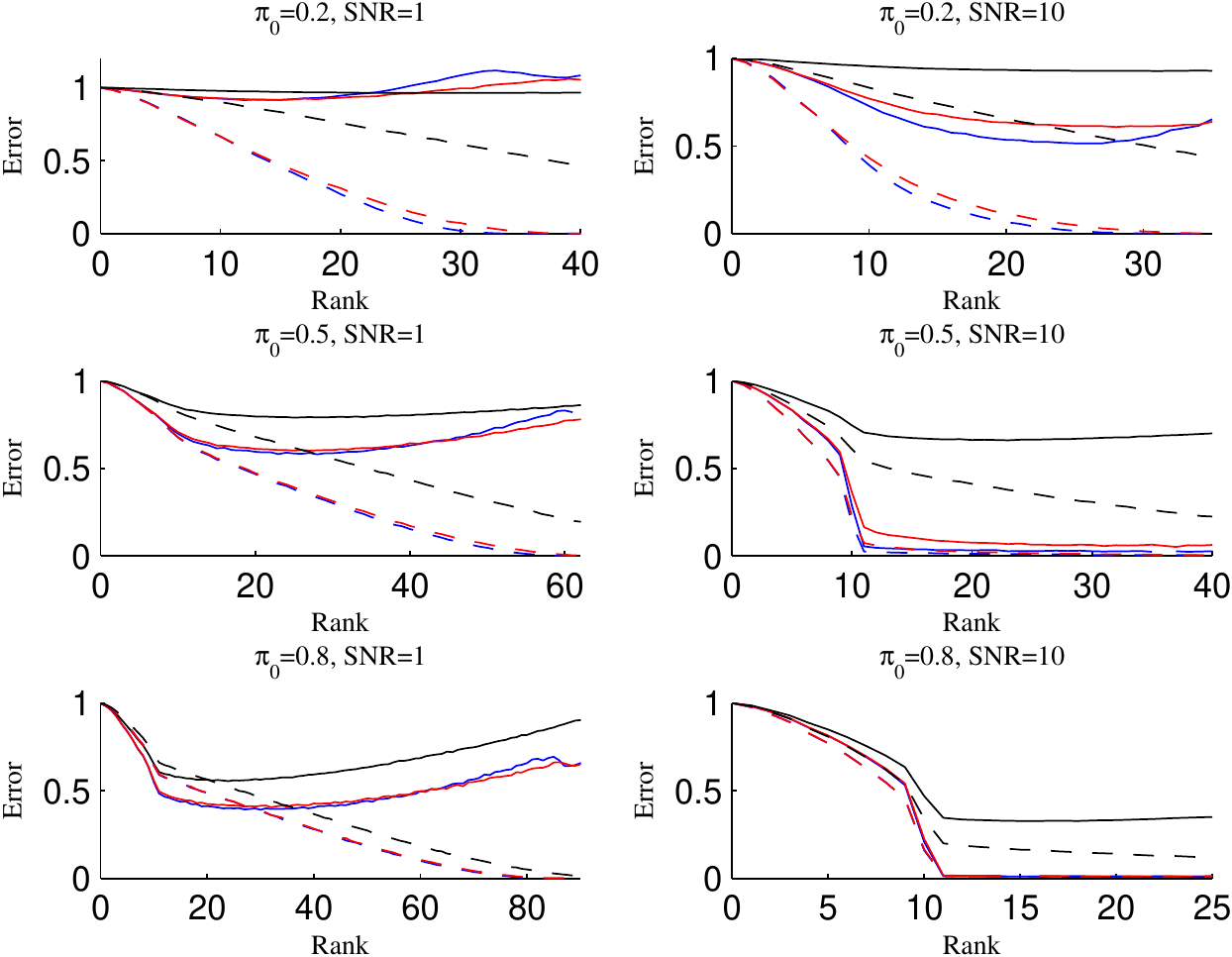}
\label{fig-simu6}
\caption{Plots of training and testing errors against the estimated rank: testing error with solid lines; training
error with dashed lines; spectrum Lasso in blue, calibrated spectrum E-net in red; modified spectrum Lasso in black;
$d_1=d_2=100$, $\rank(\Theta)=10$.}
\end{figure}

We compare the calibrated spectrum E-net (\ref{cali-elastic}) with the spectrum Lasso (\ref{matrix-lasso})
and its modification $\Thetahat^{(\text{KLT})}$ of \cite{KolTL10}.
For all methods, we compute a series of estimators with 100 different penalty levels, where the smallest penalty level corresponds to a full-rank solution and the largest penalty level corresponds to a zero solution. For the calibrated spectrum E-net, we always use $\lam_2=\lam_1\{n/(d\log d)\}^{1/4}/\Fhat $, where $\Fhat =(\sum_{i=1}^ny_i^2/\pi_0)^{1/2}$ is an estimator for $\|\Theta\|_{(F)}$.
We plot the training errors and test errors as functions of estimated ranks, where the training and test errors are defined as
\bes
\text{Training error}=\frac{\|\calP_{\Omega}(\Thetahat-Y)\|_{(F)}^2}{\|\calP_{\Omega}Y\|_{(F)}^2},\quad
\text{Test error}=\frac{\|\calP_{\Omega}^\perp(\Thetahat-\Theta)\|_{(F)}^2}{\|\calP_{\Omega}^\perp\Theta\|_{(F)}^2}.
\ees

In Figure 1, we report the estimation performance of three methods. The rank of $\Theta$ is 10 but
$\{\Theta, \Omega,\veps\}$ are regenerated in each replication.
Different noise levels and proportions of the observed entries are considered.
All the results are averaged over 50 replications.
In this experiment, the calibrated spectrum E-net and the spectrum Lasso estimator have very close testing
and training errors, and both of them significantly outperform the modified Lasso.
Figure 1 also illustrates that in most cases, the calibrated spectrum E-net and spectrum Lasso achieve
the optimal test error when the estimated rank is around the true rank.

We note \chzblue{that} 
the constrained spectrum Lasso
estimator $\Thetahat^{(\text{NW})}$ would have the same performance \chzblue{as the spectrum Lasso} 
when the constraint
$\alpha_{sp}(\Thetahat)\le\alpha^*$ is set with a sufficiently high $\alpha^*$. However, analytic properties of
the spectrum Lasso is unclear without constraint or modification.

\section{Proof of Theorem \ref{thm-oracle}}
The proof of Theorem \ref{thm-oracle} requires the following proposition that controls the approximation error
of the Taylor expansion of the nuclear norm with subdifferentiation.
The result, closely related to those in \cite{Wedin72}, is used to control the variation of the tangent space
of the spectrum E-net estimator. We omit its proof.

\begin{proposition}\label{prop-sub-diff}
Let $\Theta = UDV^\top$ be the SVD and $M$ be another matrix. Then,
\bes
0 &\le& \|M\|_{(N)} - \|\Theta\|_{(N)} - \|\calP_T^\perp M\|_{(N)} - \langle UV^\top, M-\Theta \rangle
\cr &\le& \| (\calP_TM-\Theta)VD^{-1/2}\|_{(F)}^2 + \| D^{-1/2}U^\top (\calP_TM-\Theta)\|_{(F)}^2.
\ees
\end{proposition}

\noindent
{\bf Proof of Theorem \ref{thm-oracle}.}
Define
\bes
&&\Theta^* = (\calP_T\calH\calP_T+\lam_2\calP_T)^{-1}(\calP_T\veps + \calP_T\calH\Theta - \lam_1 {U}{V}^\top),\\
&&\Thetabar = (\pi_0+\lam_2)^{-1}(\pi_0\Theta - \lam_1 {U}{V}^\top),\\
&&\Delta = \Thetatil - \Theta^*,\quad \Delta^* = \Theta^*-\Thetabar,\quad \Delta_*=\Thetatil - \Thetabar.
\ees
Since $\Thetahat = \xi\Thetatil$ and $\xi\Thetabar - \Theta = - (\lam_1/\pi_0)UV^\top$,
\bel{pf-case2}\nonumber
\|\Thetahat - \Theta\|_{(F)}&\le& \xi\|\Delta_*\|_{(F)} + \|{\xi}\Thetabar  - \Theta\|_{(F)}
\\ &=&\xi\|\Delta_*\|_{(F)} + \sqrt{r}\lam_1/\pi_0
\\ &\le& \xi\|\Delta\|_{(F)} + \xi\|\Delta^*\|_{(F)} + \sqrt{r}\lam_1/\pi_0\label{pf-case1}.
\eel
We consider two cases by comparing $\lam_2$ and $\pi_0$.

{\it Case 1: $\lam_2\le\pi_0$.}
By algebra $\xi\Delta^* = \pi_0^{-1}(\calP_T\calR+\calP_T)^{-1}\calP_T(\veps + \Deltabar)$, so that
\bel{pf-Delta-star}
\xi\|\Delta^*\|_{(F)} \le \pi_0^{-1} \|(\calP_T\calR+\calP_T)^{-1}\|_{(op)}\|\calP_T\Deltabar+\calP_T\veps\|_{(F)}
\le \sqrt{r}\lam_1/(2\pi_0). 
\eel
The last inequality above follows from the first inequalities in (\ref{event1}), (\ref{event2}) and (\ref{event3}).
It remains to bound $\|\Delta\|_{(F)}$. Let $Y  = \sum_{i=1}^n y_i E_{\omega_i}$. We write the spectrum E-net
estimator (\ref{Thetatil}) as
\bes
\Thetatil = \argmin_M\Big\{\langle \calH M,M\rangle/2 - \langle Y,M\rangle
+ \lam_1\|M\|_{(N)}+(\lam_2/2)\|M\|_{(F)}^2\Big\}.
\ees
It follows that for a certain member $\Ghat$ in the sub-differential of $\|M\|_{(N)}$ at $M=\Thetatil$,
\bes
0 = \pa L_{\lam_1,\lam_2}(\Thetatil)
= \calH \Thetatil - Y + \lam_2\Thetatil + \lam_1\Ghat
= (\calH +\lam_2)\Delta + (\calH+\lam_2)\Theta^* - Y + \lam_1\Ghat.
\ees
Let $\Rem_1 = \|\Theta^*\|_{(N)}-\langle {U}{V}^\top,\Theta^*\rangle$.
Since $\|\Theta^*\|_{(N)} - \|\Thetatil\|_{(N)} \ge - \langle \Delta, \Ghat\rangle$, we have
\bel{pf-Delta-1}
\langle (\calH +\lam_2)\Delta,\Delta\rangle
&\le& \langle \calH\Theta+\veps - (\calH+\lam_2)\Theta^*,\Delta\rangle
+\lam_1\|\Theta^*\|_{(N)} - \lam_1\|\Thetatil\|_{(N)}
\cr  &=& \langle \calH(\Theta-\Theta^*) +\veps - \lam_2\Theta^*,\Delta\rangle
+\lam_1\Rem_1 + \lam_1 \langle UV^\top, \Theta^*\rangle - \lam_1\|\Thetatil\|_{(N)}
\cr  &\le& \lam_1 \Rem_1 + \langle \veps + \calH(\Theta-\Theta^*)-\lam_2\Theta^* - \lam_1 {U}{V}^\top, \Delta\rangle
- \lam_1 \|\calP_T^\perp\Delta\|_{(N)} 
\cr  &=& \lam_1 \Rem_1 + \langle \veps+\calH(\Theta-\Theta^*),
\calP_T^\perp \Delta\rangle- \lam_1 \|\calP_T^\perp\Delta\|_{(N)}.
\eel
The second inequality in (\ref{pf-Delta-1}) is due to $\|\Thetatil \|_{(N)}
\ge \|\calP_T^\perp \Thetatil \|_{(N)} + \langle UV^\top,\Thetatil \rangle$ and $\calP_T^\perp \Thetatil=\calP_T^\perp \Delta$.
The last equality in (\ref{pf-Delta-1}) follows from the definition of $\Theta^*\in T$, since it gives
$\calP_T\veps + \calP_T\calH(\Theta-\Theta^*)-\lam_2\Theta^* - \lam_1 {U}{V}^\top=
 - (\calP_T\calH\calP_T+\lam_2\calP_T)\Theta^*+\calP_T\veps + \calP_T\calH\Theta - \lam_1 {U}{V}^\top=0$.
By the definitions of $\calQ$, $\Theta^*$ and $\Deltabar$,
$\veps+\calH(\Theta-\Theta^*)
= \calQ\veps+\calH(\Theta-\Thetabar)
- \calH(\calP_T\calH\calP_T+\lam_2\calP_T)^{-1}\calP_T\Deltabar$.
Since \chzblue{$\calP_T^\perp\calH\calP_T=\calP_T^\perp(\calH-\pi_0)\calP_T=\calP_T^\perp\calR(\pi_0+\lam_2)$}
and $(\calH-\pi_0)(\Theta-\Thetabar) = \Deltabar$, we find
\bes
&& \langle \veps+\calH(\Theta-\Theta^*), \calP_T^\perp \Delta\rangle
\cr&= & \langle \calQ\veps+(\calH-\pi_0)\{\Theta-\Thetabar
- (\calP_T\calH\calP_T+\lam_2\calP_T)^{-1}\calP_T\Deltabar\}, \calP_T^\perp \Delta\rangle
\cr&= & \langle \calQ\veps+\Deltabar - \calR(\calP_T\calR+\calP_T)^{-1}\calP_T\Deltabar, 
\calP_T^\perp \Delta\rangle. 
\ees
Thus, by the second inequalities of (\ref{event2}) and (\ref{event3}),
\bel{pf-veps}
\langle \veps+\calH(\Theta-\Theta^*), \calP_T^\perp \Delta\rangle \le \lam_1\|\calP_T^\perp \Delta\|_{(N)}.
\eel
Since $\Theta^* = \Delta^*-\Thetabar \in T$ and the singular values of $\Thetabar$ is no smaller than
$(\pi_0 s_r - \lam_1)/(\pi_0+\lam_2)\ge (s_r - \lam_1/\lam_2)/\xi \ge 4\lam_1/(\lam_2\xi)$
by the second inequality in (\ref{event1}), Proposition \ref{prop-sub-diff} and (\ref{pf-Delta-star}) imply
\bel{pf-rem}
\Rem_1 \le 2\|\Theta^* - \Thetabar\|_{(F)}^2/\{(\pi_0s_r-\lam_1)/(\pi_0+\lam_2)\} \le r(\lam_1/\pi_0)^2/(8\xi\lam_1/\lam_2).
\eel
It follows from (\ref{pf-Delta-1}), (\ref{pf-veps}) and (\ref{pf-rem}) that 
\bel{pf-Delta}
\xi^2\|\Delta\|_{(F)}^2 \le \xi^2\langle (\calH +\lam_2)\Delta,\Delta\rangle/\lam_2
\le \xi^2(\lam_1/\lam_2)\Rem_1 
\le r\lam_1^2/(4\pi_0^2).
\eel
Therefore, the error bound (\ref{oracle-bound}) follows from (\ref{pf-case1}), (\ref{pf-Delta-star}) and (\ref{pf-Delta}).

{\it Case 2: $\lam_2\ge\pi_0$.} By applying the derivation of (\ref{pf-Delta-1}) to $\Thetabar$ instead of $\Theta^*$, we find
\bes
&& \langle (\calH+\lam_2)\Delta_* , \Delta_*\rangle + \lam_1 \|\calP_T^\perp\Delta_*\|_{(N)}
\cr &\le& \lam_1\big(\|\Thetabar\|_{(N)}-\langle UV^\top,\Thetabar\rangle\big)
+ \langle \veps + \calH(\Theta-\Thetabar)-\lam_2\Thetabar - \lam_1 {U}{V}^\top, \Delta_*\rangle.
\ees
By the definitions of $\Deltabar$, $\calR$, and $\Thetabar$,
$\Deltabar 
=(\calH - \pi_0)(\Theta-\Thetabar)
=\calH(\Theta-\Thetabar)-\lam_2\Thetabar - \lam_1 {U}{V}^\top$.
This and $\|\Thetabar\|_{(N)}=\langle UV^\top,\Thetabar\rangle$ gives
\bel{pf-Delta-lowstar-1}
\langle (\calH+\lam_2)\Delta_* , \Delta_*\rangle + \lam_1 \|\calP_T^\perp\Delta_*\|_{(N)}
\le \langle \veps + \Deltabar, \Delta_*\rangle.
\eel
Since $\|\calP_T^\perp(\veps+\Deltabar)\|_{(S)}=\|\calP_T^\perp\veps\|_{(S)}\le\lam_1$
by the third inequality in (\ref{event3}), we have
\bel{pf-Delta-lowstar-2}
\langle \calP_T^\perp(\veps + \Deltabar), \Delta_*\rangle
\le \lam_1\|\calP_T^\perp\Delta_*\|_{(N)}.
\eel
It follows from (\ref{pf-Delta-lowstar-1}), (\ref{pf-Delta-lowstar-2}) and the first inequalities of
(\ref{event2}) and (\ref{event3}) that
\bes
\lam_2\|\Delta_*\|_{(F)}^2
\le \langle \calP_T(\veps + \Deltabar), \Delta_*\rangle
\le \Big\{\|\calP_T\veps\|_{(F)}+\|\calP_T\Deltabar\|_{(F)}\Big\}\|\Delta_*\|_{(F)}
\le \sqrt{r}\lam_1\|\Delta_*\|_{(F)}/2.
\ees
Thus, due to $\lam_2\ge\pi_0$,
\bel{pf-Delta-lowstar}
\xi\|\Delta_*\|_{(F)}\le (\xi/\lam_2)\sqrt{r}\lam_1/2\le \sqrt{r}\lam_1/\pi_0.
\eel
Therefore, the error bound (\ref{oracle-bound}) follows from (\ref{pf-case2}) and (\ref{pf-Delta-lowstar}).
$\hfill\square$

\subsubsection*{Acknowledgments}
This research is partially supported by the NSF Grants DMS 0906420, DMS-11-06753 
\chzblue{and DMS-12-09014}, 
and NSA Grant H98230-11-1-0205.

\newpage
\bibliographystyle{abbrv}
\bibliography{MatrixCEN}

\end{document}

%% file: 04-pre-am-short.tex
\usepackage{amsmath}
\usepackage{latexsym}
\usepackage{mathrsfs} 
\usepackage{amsfonts}
\usepackage{amssymb}

\newcommand{\bel}{\begin{eqnarray}\label}
\newcommand{\eel}{\end{eqnarray}}
\newcommand{\bes}{\begin{eqnarray*}}
\newcommand{\ees}{\end{eqnarray*}}

\def\benu{\begin{enumerate}}
\def\eenu{\end{enumerate}}

\def\argmin{\mathop{\rm arg\, min}}

\def\real{{\mathbb{R}}}

\def\complex{\mathop{{\rm I}\kern-.58em\hbox{\rm C}}\nolimits}
\def\pa{\partial}

\def\diag{\hbox{diag}}

\def\rank{\hbox{rank}}
\def\Rem{\hbox{Rem}}

\def\trace{\hbox{trace}}

\newcommand{\R}{{\mathbb{R}}}

\def\Fhat{\widehat{F}}

\def\Ghat{\widehat{G}}

\def\calH{{\cal H}}

\def\calI{{\cal I}}

\def\calP{{\cal P}}

\def\calQ{{\cal Q}}

\def\calR{{\cal R}}

\def\Ybar{{\overline Y}}

\def\Deltabar{{\overline \Delta}}

\def\eps{\epsilon}\def\veps{\varepsilon}

\def\Thetahat{\widehat{\Theta}}
\def\Thetatil{{\widetilde \Theta}}\def\Thetabar{{\overline \Theta}}

\def\lam{\lambda}